\documentclass[11pt]{article}
\usepackage{amsfonts}
\usepackage{latexsym,amssymb,amsmath,graphics,cite,arydshln}
\usepackage{tikz}
\usepackage{fancyhdr}
\usepackage{lipsum}
\usepackage{lineno}
\usepackage{color}

\usepackage{algorithm} %format of the algorithm
\usepackage{algorithmic} %format of the algorithm
\usepackage{multirow} %multirow for format of table
\usepackage{amsmath}
\usepackage{xcolor}

\begin{document}
\newcommand{\qed}{\hphantom{.} \hfill $\Box$\medbreak}
\newcommand{\proof}{\noindent{\bf Proof \ }}

\newtheorem{theorem}{Theorem}[section]
\newtheorem{lemma}[theorem]{Lemma}
\newtheorem{corollary}[theorem]{Corollary}
\newtheorem{remark}[theorem]{Remark}
\newtheorem{example}[theorem]{Example}
\newtheorem{definition}[theorem]{Definition}
\newtheorem{fact}[theorem]{Fact}
\newtheorem{claim}[theorem]{Claim}
\newtheorem{proposition}[theorem]{Proposition}
\newtheorem{note}[theorem]{Note}

\begin{center}
{\Large\bf Several problems on reduced spherical polygons of thickness less than $\pi/2$\footnote{Supported by NSFC under Grant $11971053$.}
}

\vskip12pt
Cen Liu and Yanxun Chang\\[2ex] {\footnotesize Department of Mathematics, Beijing Jiaotong University, Beijing 100044, China}

{\footnotesize  cenliumath@foxmail.com, yxchang@bjtu.edu.cn}
\vskip12pt
\end{center}

\vskip12pt

\noindent {\bf Abstract:}
The present paper aims to solve some problems proposed by Lassak about the reduced spherical polygons. The main result is to show that the regular spherical $n$-gon has the minimal perimeter among all reduced spherical polygons of fixed thickness less than $\pi/2$ and with at most $n$ vertices.
In addition, we determine the maximal diameter of every reduced spherical polygons with a fixed thickness less than $\pi/2$. We also find the smallest spherical radius that contains every reduced spherical polygons with a fixed thickness less than $\pi/2$.
\vskip12pt


%%%%%%%%%%%%%%%%%%%%%%%%%%%%%%%%%%%%%%%%%%%%%%%%%%%%%%%%%%%%%%%%%%%%%%%%%

\section{Introduction}
Reduced convex bodies have been investigated in Euclidean space, normed space and spherical space, respectively, one can refer to \cite{convexbody,survey2,survey3,survey4}.
The present paper focus on the reduced polygons on the sphere.
Particularly, the investigation of reduced spherically convex polygon was first introduced by Lassak in \cite{polygon}.

One can check that, the last paragraph of the article \cite{polygon} establishes four conjectures about the perimeter and area of reduced spherically convex polygons with a fixed thickness less than $\pi/2$.
In particular, the four conjectures have been solved in \cite{CL} with a given thickness $\pi/2$.
Furthermore, two conjectures about the area of reduced spherical polygon of thickness less than $\pi/2$ have been solved in \cite{CenL}. The main result of this paper confirms one of the remaining conjectures which shows that the regular spherical $n$-gon has the minimal perimeter among all reduced spherical $n$-gons of a fixed thickness less than $\pi/2$.
In addition, there are still some other problems about the reduced spherical polygons proposed by Lassak in \cite{polygon} and \cite{survey4}, respectively. We consider these problems in the following.

\begin{itemize}
\item[$\mathrm{P1}$.] Is it true that the regular spherical $n$-gon has minimal perimeter among all reduced spherical polygons of the same thickness less than $\pi/2$ and with at most $n$ vertices?
\item[$\mathrm{P2}$.] Is it true that the equality in Theorem $4.2$ \cite{polygon} holds only for regular triangle?
\item[$\mathrm{P3}$.] What is the smallest radius of a disk which contains every reduced polygon of a given thickness on $S^{2}$?
\end{itemize}

In Section $2$, we review some results from the literature and several useful lemmas are established.
Section $3$ confirms the Problem $\mathrm{P1}$. Actually, we use a similar method to solve this problem as the proof process in \cite{CenL}.
Although the answer to Problem $\mathrm{P2}$ in Section $4$ is negative. We determine the supremum of the diameter of any reduced spherically convex polygons and confirm $\mathrm{P2}$ in our investigation.
In Section $5$, based on the consideration of $\mathrm{P2}$,
we get the smallest radius of a disk that contains every reduced spherical polygons of a fixed thickness less than $\pi/2$, which answers the Problem $\mathrm{P3}$.

\section{Preliminaries}
For the notions of the reduced convex polygons on $S^{2}$, one can see \cite{CL,polygon,survey4,CenL}.
Let $S^{2}$ be the unit sphere of $E^{3}$. The intersection of $S^{2}$ with any two-dimensional subspace of $E^{3}$ is called a {\it great circle}.
 A pair of {\it antipodes} are the intersection of $S^{2}$ with any one-dimensional subspace of $E^{3}$.
If $a,b\in S^{2}$ are different points and not antipodes, then there is exactly one great circle containing them; Denote by {\it arc ab}, shortly $ab$, the shorter part of the great circle containing them. The {\it distance} $|ab|$ is the length of $ab$.

 Let $C$ be a spherical set, for any two points $a,b\in C$, if we have $ab\subseteq C$ and $|ab|<\pi$, then $C$ is a {\it convex} set.
The {\it convex body} is a closed convex set with non-empty interior.

The {\it spherical disk} of radius $r\in(0,\pi/2]$ and {\it center} $k\in S^{2}$ is the set of points having distance at most $r$ from $k$;
the {\it spherical circle} is the boundary of the spherical disk.
Spherical disks of radius $\pi/2$ are called {\it hemispheres}.
If $G$ and $H$ are different hemispheres and their centers are not antipodes,
then $L=G\cap H$ is called a {\it lune} of $S^{2}$.
The parts of $\mathrm{bd}(G)$ and $\mathrm{bd}(H)$ contained in $G\cap H$ are
denoted by $G/H$ and $H/G$, respectively.
 We define the {\it thickness} $\Delta(L)$ of the lune $L=G\cap H $ as the distance between the centers of $G/H$ and $H/G$. The thickness of a spherically convex polygon $C$ is the minimal thickness of a lune which contains this polygon, and we denote it by $\Delta(C)$.
A convex body $C\subset S^{2}$ is said to be {\it reduced}
if $\Delta(R)<\Delta(C)$ for each convex body $R$ being a proper subset of $C$.

Some useful definitions given in \cite{polygon} are established here.
If $C$ is a subset of a convex set of $S^{2}$, then the intersection of all convex sets containing $C$ is called a {\it convex hull} of $C$.
The convex hull of $k\geq 3$ points on $S^{2}$
such that each of them does not belong to the convex hull of the remaining
points is called a {\it spherically convex} $k$-gon.
If $V$ is a spherically convex $k$-gon, we denote by $v_{1},\ldots ,v_{k}$ the vertices of $V$ in the counterclockwise order. We define the counterclockwise direction as the positive orientation.
A spherically convex polygon with sides of equal length and interior angles of equal measure is called a {\it regular spherical polygon}.

Let $p$ be a point in a hemisphere different from its center and let $l$ be the great circle bounding this hemisphere. The {\it projection} of $p$ on $l$ is the point $t$ such that $|pt|=\min \{|pc|:c\in l\}$.
For a spherically convex odd-gon $V=v_{1}v_{2}\cdots v_{n}$, by the {\it opposite side to the vertex $v_{i}$} we mean the side $v_{i+(n-1)/2}v_{i+(n+1)/2}$, where the indices are taken modulo $n$.

A few formulas of spherical trigonometry in \cite{DA} are useful for our research.
For a right spherical triangle with hypotenuse $c$ and legs $a,b$,
denote by $A,B$ and $C$ the corresponding angles of edges $a,b$ and $c$, respectively. Then we have
\begin{equation}\label{eq1}
\cos A=\tan b\cot c,
\end{equation}
\begin{equation}\label{eq-1}
\cos B=\tan a\cot c,
\end{equation}
\begin{equation}\label{eq-3}
\sin b=\sin c\sin B,
\end{equation}
\begin{equation}\label{eq2}
\cos c=\cos a\cos b,
\end{equation}
\begin{equation}\label{eq-2}
\cos c=\cot A\cot B,
\end{equation}
\begin{equation}\label{eq3}
\cos B=\cos b\sin A.
\end{equation}
For any spherical triangle, by the sine theorem on sphere, we have
\begin{equation}\label{eq-4}
\frac{\sin A}{\sin a}=\frac{\sin B}{\sin b}=\frac{\sin C}{\sin c}.
\end{equation}
Our paper is based on the following results given in \cite{polygon} and \cite{CenL}, respectively.

\begin{lemma}{\rm\cite[Theorem 3.2]{polygon}}\label{3.2}
Every reduced spherical polygon is an odd-gon of thickness at most $\frac{\pi}{2}$.
A spherically convex odd-gon $V$ with $\Delta(V)<\frac{\pi}{2}$ is reduced if and only if the projection of every its vertices on the great circle containing the opposite side belongs to the relative interior of this side and the distance of this vertex from this side is $\Delta(V)$.
\end{lemma}

\begin{lemma}{\rm\cite[Corollary 3.3]{polygon}}\label{col3.3}
Every spherical regular odd-gon of thickness at most $\frac{\pi}{2}$ is reduced.
\end{lemma}

\begin{lemma}{\rm\cite[Corollary 3.6]{polygon}}\label{3.6}
For every reduced odd-gon $V=v_{1}v_{2}\cdots v_{n}$ with $\Delta(V)<\frac{\pi}{2}$, we have
$|v_{i}t_{i+(n+1)/2}|=|t_{i}v_{i+(n+1)/2}|$, for $i=1,2,\ldots,n$, where $t_{i}$ denotes the projection of $v_{i}$ on the opposite side.
\end{lemma}

\begin{lemma}{\rm\cite[Corollary 3.9]{polygon}}\label{col3.9}
If $V=v_{1}v_{2}\cdots v_{n}$ is a reduced spherical polygon with $\Delta(V)<\frac{\pi}{2}$,
then $\beta_{i}\leq\alpha_{i}$ for every $i\in\{1,\ldots,n\}$.
\end{lemma}

\begin{figure}[h]
\centering
  \includegraphics[width=7 cm]{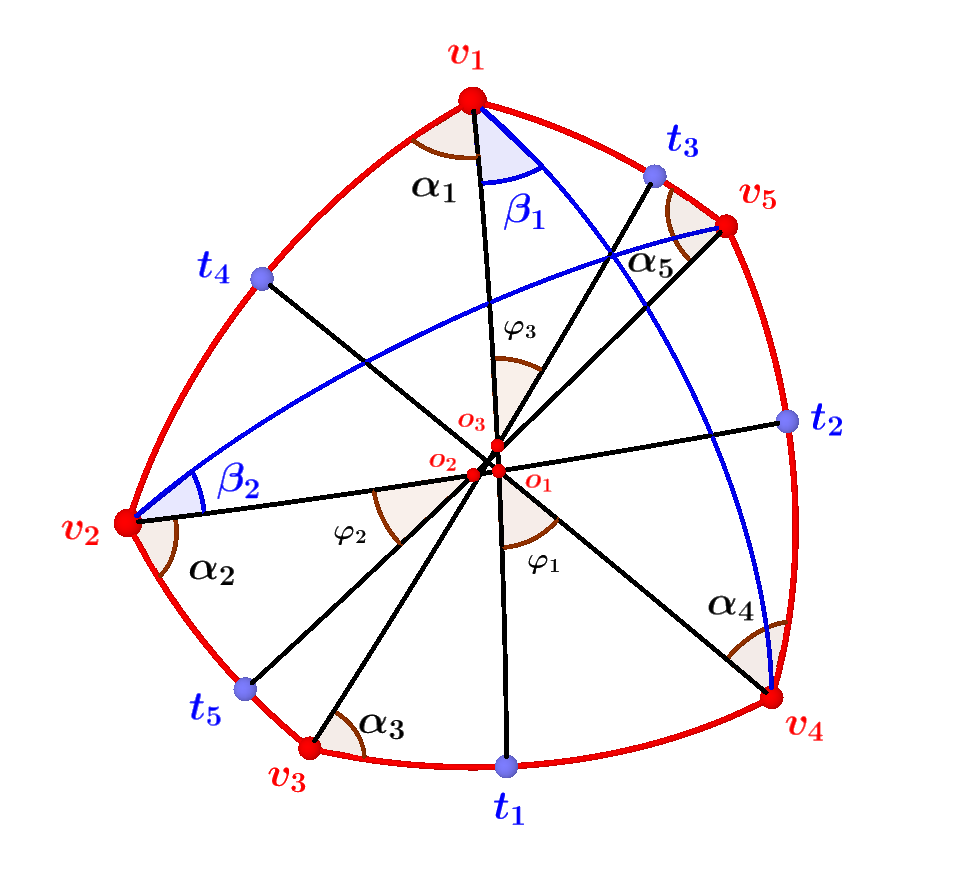}\\
\caption{Related notations in a reduced spherical pentagon}
  \label{example}
\end{figure}

We use the similar method as in \cite{CenL} to prove the Problem $\mathrm{P1}$, then the notations given in \cite[Section $3$]{CenL} are applicable and we rewrite it here.
In a reduced spherical polygon $V=v_{1}v_{2}\cdots v_{n}$, let $t_{i}$ be the projection of $v_{i}$ on the opposite side $v_{i+(n-1)/2}v_{i+(n+1)/2}$.
 Let $o_{i}$ be the intersection point of $v_{i}t_{i}$ and $v_{i+(n+1)/2}t_{i+(n+1)/2}$.
Put $\alpha_{i}=\angle v_{i+1}v_{i}t_{i}$, $\beta_{i}=\angle t_{i}v_{i}v_{i+(n+1)/2}$,
 and $\varphi_{i}=\angle v_{i}o_{i}t_{i+(n+1)/2}=\angle t_{i}o_{i}v_{i+(n+1)/2}$, where $i\in\{1,2,\ldots,n\}$.
Fig.~\ref{example} presents some notations in a reduced spherical pentagon.

\begin{lemma}{\rm\cite[Lemma 3.2]{CenL}}\label{cog-tri}
For every reduced spherical odd-gon $V=v_{1}v_{2}\cdots v_{n}$ with $\Delta(V)<\frac{\pi}{2}$, the spherical triangles $v_{i}o_{i}t_{i+(n+1)/2}$ and $v_{i+(n+1)/2}o_{i}t_{i}$ are congruent, where $i\in\{1,2,\ldots,n\}$.
\end{lemma}

\begin{lemma}{\rm\cite[Fact 3.4]{CenL}}\label{fact-2020}
For every reduced spherical polygon $V=v_{1}v_{2}\cdots v_{n}$ with a thickness $\Delta(V)<\frac{\pi}{2}$, we have $0<\varphi_{i}<\frac{\pi}{2}$, where $i\in\{1,\ldots,n\}$.
\end{lemma}

\begin{lemma}{\rm\cite[Lemma 3.5]{CenL}}\label{lem3.5}
For any reduced spherical polygon $V=v_{1}v_{2}\cdots v_{n}$ of thickness less than $\frac{\pi}{2}$, we have $\sum_{i=1}^n\varphi_{i}\geq\pi$. Moreover,
\begin{itemize}
\item[$(1)$]  if $V$ is a non-regular reduced spherical polygon, then $\sum_{i=1}^n\varphi_{i}\geq\pi$;
\item[$(2)$]  if $V$ is a regular spherical polygon, then $\sum_{i=1}^n\varphi_{i}=\pi$ and $\varphi_{i}=\frac{\pi}{n}$, where $i\in\{1,2,\ldots,n\}$.
\end{itemize}
\end{lemma}

In the following, we investigate the length of side, interior angle and circumradius of regular spherical triangle with a given thickness. Since different formulas of spherical trigonometry derive different expressions, here we establish one of them.

\begin{lemma}\label{gamma-regular triangle}
Let $V=v_{1}v_{2}v_{3}$ be a regular spherical triangle with thickness less than  $\frac{\pi}{2}$. Then we have
 \begin{itemize}
\item[$(1)$] the interior angle of $V$ is $2\gamma$;
\item[$(2)$] the length of one side of $V$ is $2\arccos(\frac{1}{2\sin\gamma})$;
\item[$(3)$] the circumradius of $V$ is $\arcsin(\frac{2}{\sqrt{3}}\sqrt{1-\frac{1}{4\sin^{2}\gamma}})$,
\end{itemize}
where $\gamma=\arcsin\frac{-\cos\Delta(V)+\sqrt{\cos^{2}\Delta(V)+8}}{4}$.
\end{lemma}
\proof
Lemma \ref{col3.3} permits that $V$ is reduced. Hence we apply the above notations to $V$.

From the knowledge of regular spherical trigonometry, we know that $\alpha_{i}=\beta_{i}$ and $t_{i}$ is the midpoint of $v_{i+1}v_{i+2}$, where $i=1,2,3$ and the indices are taken modulo $3$.
Consider the right spherical triangle $v_{i}t_{i}v_{i+2}$, then we have $|v_{i}t_{i}|=\Delta(V)$.
Let $|t_{i}v_{i+2}|=a$, thus $|v_{i}v_{i+2}|=2a$.
Denote by $\gamma$ the angle $\angle t_{i}v_{i}v_{i+2}$.
By Equation \eqref{eq2}, we get that
\begin{equation}\label{regular-qru}
\cos 2a=\cos\Delta(V)\cos a,
\end{equation}
and thus $\cos a=\frac{\cos\Delta(V)+\sqrt{\cos^{2}\Delta(V)+8}}{4}$.
By Equation \eqref{eq-4}, we get that
$$
\frac{\sin\gamma}{\sin a}=\frac{\sin\frac{\pi}{2}}{\sin2a},
$$
from this and \eqref{regular-qru}, we obtain that $$\gamma=\arcsin\frac{-\cos\Delta(V)+\sqrt{\cos^{2}\Delta(V)+8}}{4}.$$
The interior angle of $V$ is as desired.
Since $\Delta(V)\in(0,\frac{\pi}{2})$, it follows that $\gamma\in(\frac{\pi}{6},\frac{\pi}{4})$ and $2\gamma\in(\frac{\pi}{3},\frac{\pi}{2})$.

 Consider the right spherical triangle $v_{1}t_{2}o$ in $V$, where $o$ is the center point of the circumcircle of $V$. By using \eqref{eq3}, we obtain that the length of $v_{1}t_{2}$ is $\arccos(\frac{1}{2\sin\gamma})$. Hence the length of one side of $V$ is $2\arccos(\frac{1}{2\sin\gamma})$.
From Equation \eqref{eq-3}, it follows that the circumradius of $V$ is $\arcsin(\frac{2}{\sqrt{3}}\sqrt{1-\frac{1}{4\sin^{2}\gamma}})$.
\hfill\qed

In fact, $\gamma$ is half of the interior angle of a regular triangle of thickness $\omega\in(0,\frac{\pi}{2})$. For convenient, we denote the number $\arcsin\frac{-\cos\omega+\sqrt{\cos^{2}\omega+8}}{4}$ by $\gamma$ in the whole research, where $\gamma\in(\frac{\pi}{6},\frac{\pi}{4})$.

By Lemma \ref{col3.9}, for each reduced spherical $n$-gon $V$, it is true that $\beta_{i}\leq\alpha_{i}$, (if $V$ is a reduced polygon in the plane, the relation is $\beta_{i}\leq\frac{\pi}{6}\leq\alpha_{i}$, one can refer \cite[Theorem $8$]{1990}). Here we intend to find a similar relation between $\beta_{i}$ and $\alpha_{i}$ as in the plane.

In a spherically convex polygon $V$, for arbitrary two points $x,y\in\mathrm{bd}(V)$, denote by $\widehat{xy}$ the arcs of $\mathrm{bd}(V)$ from $x$ to $y$ in a positive orientation. By $|\widehat{xy}|$ we mean the length of $\widehat{xy}$.

\begin{lemma}\label{relation}
For each reduced spherical polygon $V=v_{1}v_{2}\cdots v_{n}$ of thickness less than $\frac{\pi}{2}$, we have $\beta_{i}\leq \gamma\leq\alpha_{i}$, where $i\in\{1,\ldots,n\}$.
Especially, $\beta_{i}$ and $\alpha_{i}$ attain $\gamma$ at the same time and
only when $V$ is a regular triangle, we can get that $\beta_{i}=\gamma=\alpha_{i}$.
\end{lemma}
\proof
Consider the right spherical triangle $v_{i}t_{i}v_{i+(n+1)/2}$ in $V$.
Lemma \ref{3.2} implies that $|v_{i}t_{i}|=\Delta(V)$.
By Lemma \ref{cog-tri}, we have $\angle v_{i}v_{i+(n+1)/2}t_{i}=\alpha_{i}+\beta_{i}$.
By Equation \eqref{eq-4}, we have
$$
\frac{\sin\frac{\pi}{2}}{\sin|v_{i}v_{i+(n+1)/2}|}=\frac{\sin(\alpha_{i}+\beta_{i})}{\sin\Delta(V)},
$$
then $\Delta(V)\leq\alpha_{i}+\beta_{i}\leq\frac{\pi}{2}$.
By Equation \eqref{eq3}, we have
\begin{equation}\label{eq4}
\cos(\alpha_{i}+\beta_{i})=\cos\Delta(V)\sin\beta_{i}.
\end{equation}
Lemma \ref{col3.9} permits that $\cos(\alpha_{i}+\beta_{i})\leq\cos2\beta_{i}$,
then \eqref{eq4} satisfies
\begin{equation}\label{eq5}
\cos\Delta(V)\sin\beta_{i}\leq\cos2\beta_{i}.
\end{equation}
As a result, we get $\sin\beta_{i}\leq\sin\gamma$ and thus $\beta_{i}\leq\gamma$.

Lemma \ref{col3.9} also permits that $\cos(\alpha_{i}+\beta_{i})\geq\cos2\alpha_{i}$, then \eqref{eq4} satisfies
\begin{equation}\label{eq6}
\cos\Delta(V)\sin\alpha_{i}\geq\cos2\alpha_{i}.
\end{equation}
Therefore, we get the inequality $\sin\alpha_{i}\geq\sin\gamma$ and thus $\alpha_{i}\geq\gamma$.

Consequently, the inequality $\beta_{i}\leq \gamma\leq\alpha_{i}$ holds.

Claim that $\beta_{i}$ and $\alpha_{i}$ attain $\gamma$ at the same time.
For some $i\in\{1,\ldots,n\}$, if $\beta_{i}=\gamma$, we can check that \eqref{eq5} becomes the following equality
\begin{equation*}
\cos\Delta(V)\sin\gamma=\cos2\gamma,
\end{equation*}
from this, Equation \eqref{eq4} becomes
\begin{equation*}
\cos(\alpha_{i}+\gamma)=\cos2\gamma.
\end{equation*}
Hence we obtain $\alpha_{i}=\gamma=\beta_{i}$.
For the case when $\alpha_{i}=\gamma$, we can also get the same conclusion.

\begin{figure}[htb]
\centering
  \includegraphics[width=6.5 cm]{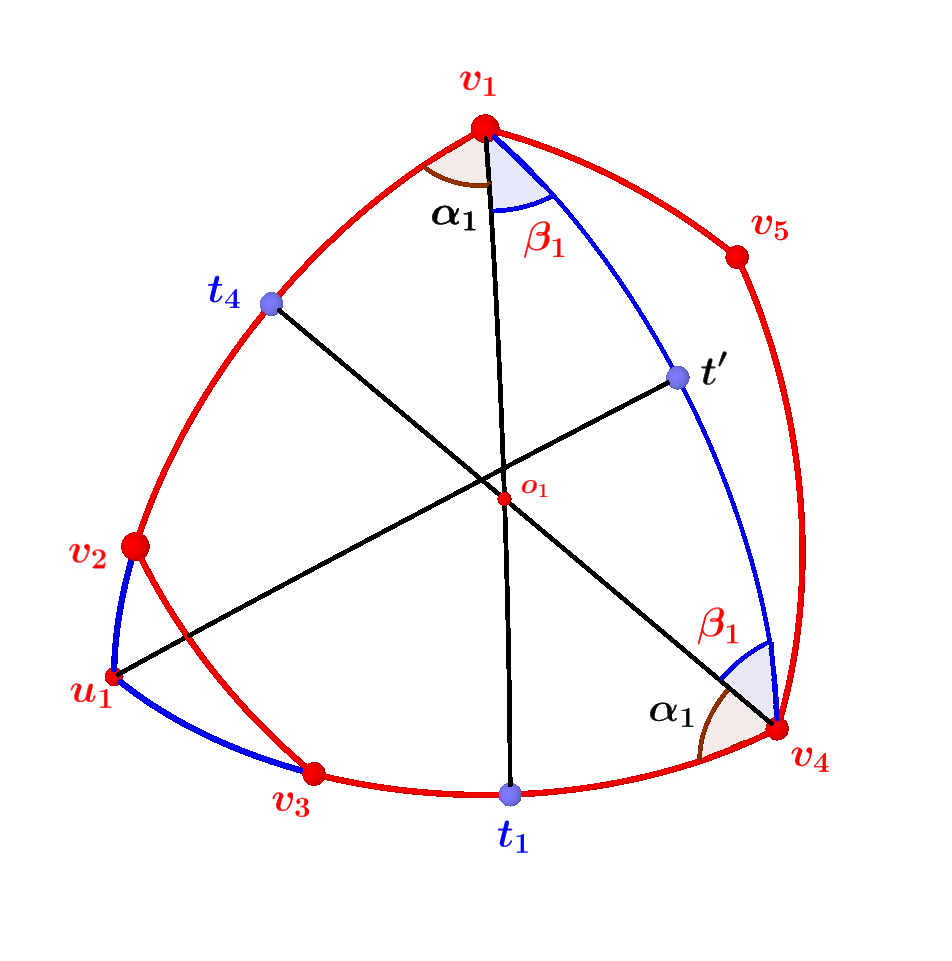}\\
\caption{Some proof details shown in a reduced spherical pentagon}
  \label{proof details}
\end{figure}

With the assumption that $\alpha_{i}=\gamma=\beta_{i}$, we will prove that $V$ degenerates into a regular spherical triangle.
Denote by $u_{i}$ the intersection point of the arcs which contains $v_{i}v_{i+1}$ and  $v_{i+(n-1)/2}v_{i+(n+1)/2}$, respectively. Fig.~\ref{proof details} illustrates the corresponding notations on a reduced spherical pentagon.

 Claim that $v_{i}u_{i}v_{i+(n+1)/2}$ is a regular spherical triangle. By the definition of regular triangle, we only need to show that the lengths of the three sides of $v_{i}u_{i}v_{i+(n+1)/2}$ are the same with each other.
 Since $\angle u_{i}v_{i}v_{i+(n+1)/2}=\angle v_{i}v_{i+(n+1)/2}u_{i}=\alpha_{i}+\beta_{i}=2\gamma$, it follows that $|u_{i}v_{i}|=|u_{i}v_{i+(n+1)/2}|$. Let $c=|u_{i}v_{i}|$ and $c'=|v_{i}v_{i+(n+1)/2}|$. Then we intend to show $c=c'$.

Let $t'$ be the projection of $u_{i}$ on the arc $v_{i}v_{i+(n+1)/2}$.
Since $|u_{i}v_{i}|=|u_{i}v_{i+(n+1)/2}|$, $t'$ is the mid-point of $v_{i}v_{i+(n+1)/2}$ and $u_{i}v_{i}t'$ is a right spherical triangle.
Then by Equation \eqref{eq-1}, we obtain $\cos 2\gamma=\tan\frac{c'}{2}\cot c$, and thus
\begin{equation}\label{eq7}
\tan c=\frac{\tan\frac{c'}{2}}{\cos 2\gamma}=\frac{\sin c'}{(1+\cos c')\cos 2\gamma}.
\end{equation}
In the right spherical triangle $v_{i}t_{i}v_{i+(n+1)/2}$. Equation \eqref{eq-2} implies that $\cos c'=\cot\gamma\cot2\gamma$, thus we get
\begin{equation*}
\frac{\cos c'}{1+\cos c'}=\frac{1}{1+\tan\gamma\tan2\gamma}
=\cos2\gamma.
\end{equation*}
Consequently, we obtain $(1+\cos c')\cos 2\gamma=\cos c'$.
From this and Equation \eqref{eq7}, it follows that $\tan c=\tan c'$. So $c=c'$, and hence $v_{i}u_{i}v_{i+(n+1)/2}$ is a regular triangle.

The following discussion is similar as that in \cite[Theorem $8$]{1990}.
From the above analysis, we get
\begin{equation}\label{eq8}
|\widehat{t_{i+(n+1)/2}t_{i}}|\leq|t_{i+(n+1)/2}u_{i}|+|u_{i}t_{i}|=|v_{i}v_{i+(n+1)/2}|\leq|\widehat{v_{i+(n+1)/2}v_{i}}|.
\end{equation}
By Lemma \ref{3.6}, we note that $$|\widehat{t_{i+(n+1)/2}t_{i}}|=|\widehat{v_{i+(n+1)/2}v_{i}}|.$$
Hence the inequality \eqref{eq8} becomes $$|\widehat{t_{i+(n+1)/2}t_{i}}|=|t_{i+(n+1)/2}u_{i}|+|u_{i}t_{i}|=|v_{i}v_{i+(n+1)/2}|=|\widehat{v_{i+(n+1)/2}v_{i}}|.$$
This means that \eqref{eq8} attains the equation only when $V$ degenerates into a regular spherical triangle. This completes the proof. \qed

The above lemma plays an important role in our consideration of the maximal diameter of reduced spherical polygons in Section $4.2$. The following discussions are very useful for Section $3$.

\begin{lemma}\label{dec-func}
Let $f_{1}(x)=\arccos\frac{1+\lambda x}{\sqrt{1+\lambda^{2}}}$ and $f_{2}(x)=\arccos\frac{x(1+\lambda x)}{\lambda-x}$. Then $\frac{f_{1}(x)}{f_{2}(x)}$ is a decreasing function of $x$, where $\lambda\in(0,+\infty)$ and $x\in(0,\frac{-1+\sqrt{1+\lambda^{2}}}{\lambda})$.
\end{lemma}
\proof We use a series of elementary calculation methods to solve this problem.
Let $f(x)=\frac{f_{1}(x)}{f_{2}(x)}$, then the first derivative of $f(x)$ is
$$
f'(x)=-\frac{\lambda \cdot h(x)}{(\lambda-x)\sqrt{-\lambda(1+x^{2})(2x+\lambda x^{2}-\lambda)}f_{2}^{2}(x)},
$$
where
$$h(x)=(\lambda-x)\sqrt{1+x^{2}}f_{2}(x)+(x^{2}-2\lambda x-1)f_{1}(x).$$
Here we have $x<\frac{-1+\sqrt{1+\lambda^{2}}}{\lambda}<\lambda$, hence $\lambda-x>0$. If $h(x)>0$ in the domain $x\in(0,\frac{-1+\sqrt{1+\lambda^{2}}}{\lambda})$, then we confirm the result $f'(x)<0$.

By a simple calculation, we get that the first derivative of $h(x)$ is
$$
h'(x)=\frac{(-2x^{2}+\lambda x-1)
f_{2}(x)}{\sqrt{1+x^{2}}}+2(x-\lambda)f_{1}(x),
$$
and the second derivative of $h(x)$ is
$$
h''(x)=\frac{\lambda(\lambda x^{3}+3x^{2}-3\lambda x+2\lambda^{2}+1)}
{(1+x^{2})(\lambda-x)\sqrt{-\lambda(2x+\lambda x^{2}-\lambda)}}
+\frac{\lambda-\frac{3x}{2}-\frac{x^{3}}{2}}{(1+x^{2})^{3/2}}f_{2}(x)+t(x),
$$
where $t(x)=-\frac{3x}{2\sqrt{1+x^{2}}}f_{2}(x)+2f_{1}(x)$.
According to the root discrimination of quadratic function, we find that
\begin{equation}\label{ineq3.8-1}
3x^{2}-3\lambda x+2\lambda^{2}+1>0.
\end{equation}
Since $\sqrt{-\lambda(2x+\lambda x^{2}-\lambda)}>0$ and $0<x<\frac{-1+\sqrt{1+\lambda^{2}}}{\lambda}<1$, it follows that
\begin{equation}\label{ineq3.8-2}
\lambda>2x+\lambda x^{2}> \frac{3x}{2}+\frac{x^{3}}{2}>0.
\end{equation}
From the inequalities \eqref{ineq3.8-1}, \eqref{ineq3.8-2}, $\lambda-x>0$ and if $t(x)>0$, then we obtain $h''(x)>0$.

By a simple calculation, we check that the first derivative of $t(x)$ is
$$
t'(x)=\frac{\lambda(x^{3}+2\lambda x^{2}+7x-4\lambda)}{2(1+x^{2})(\lambda-x)\sqrt{-\lambda(\lambda x^{2}+2x-\lambda)}}
-\frac{3}{2(1+x^{2})^{3/2}}f_{2}(x).
$$
From $\sqrt{-\lambda(\lambda x^{2}+2x-\lambda)}>0$ and $0<x<\frac{-1+\sqrt{1+\lambda^{2}}}{\lambda}<1$, we get
$$
4\lambda>4\lambda x^{2}+8x>x^{3}+2\lambda x^{2}+7x.
$$
Then $t'(x)<0$, thus we obtain $t(x)>t(\frac{-1+\sqrt{1+\lambda^{2}}}{\lambda})=0$.

Consequently, the above analysis show us $h''(x)>0$ in the domain $x\in(0,\frac{-1+\sqrt{1+\lambda^{2}}}{\lambda})$.
From $h''(x)>0$, we find that $h'(x)<h'(\frac{-1+\sqrt{1+\lambda^{2}}}{\lambda})=0$, hence $h'(x)<0$ and thus $h(x)>h(\frac{-1+\sqrt{1+\lambda^{2}}}{\lambda})=0$.
By $h(x)>0$ and $\lambda-x>0$, they confirm $f'(x)<0$, which completes the proof.
\hfill\qed

\begin{lemma}\label{F(z)}
 Let $F(x)=\arccos\frac{1+\lambda g(x)}{\sqrt{1+\lambda^{2}}}$,
where $g(x)=\frac{-(1+\cos x)+\sqrt{(1+\cos x)^{2}+4\lambda^{2}\cos x}}{2\lambda}$,
$\lambda\in(0,+\infty)$ and $x\in(0,\frac{\pi}{2})$. Then $F'(x)>0$ and $F''(x)>0$.
\end{lemma}
\proof
For convenience, set $r(x)=\sqrt{(1+\cos x)^{2}+4\lambda^{2}\cos x}$.
We find the first derivative of $F(x)$ is
\begin{equation*}
F'(x)=\frac{\cos\frac{x}{2}\sqrt{1+2\lambda^{2}+\cos x-r(x)}}
{r(x)}.
\end{equation*}
Therefore, we obtain $F'(x)>0$.  The second derivative of $F(x)$ is
\begin{equation*}
F''(x)=\frac{\lambda^{2}\sin\frac{x}{2}
(3+4\lambda^{2}+\cos x-(1-\cos x)(r(x)-\cos x))}
{r^{3}(x)\sqrt{1+2\lambda^{2}+\cos x-r(x)}}.
\end{equation*}
One can check that in the domain $x\in(0,\frac{\pi}{2})$, the following inequality holds
$$3+4\lambda^{2}+\cos x-(1-\cos x)(r(x)-\cos x)>0.$$
 Hence $F''(x)>0$ and thus $F(x)$ is a convex function of $x$.\hfill\qed

\section{The perimeter of reduced spherical polygons }
In this section, we consider the Problem $\mathrm{P1}$ proposed by Lassak in \cite{polygon}.
Let $V$ be a spherically convex $n$-gon, for simplicity, we use $\omega$ to replace $\Delta(V)$ in this section, where $\omega\in (0,\pi/2)$.
Let $\lambda=\tan\omega$, then $\lambda\in(0,+\infty)$.
Denote by $\mathrm{perim}(V)$ the perimeter of $V$.
Let us define several functions which are needed in the following theorems.
Set $$f_{1}(x)=\arccos\frac{1+\lambda x}{\sqrt{1+\lambda^{2}}}$$ and $$f_{2}(x)=\arccos\frac{x(1+\lambda x)}{\lambda-x},$$
where $x\in(0,\frac{-1+\sqrt{1+\lambda^{2}}}{\lambda})$.
Set
$$g(\varphi)=\frac{-(1+\cos\varphi)+\sqrt{(1+\cos\varphi)^{2}+4\lambda^{2}\cos \varphi}}{2\lambda},$$
where $\varphi\in(0,\frac{\pi}{2})$.
And thus $g(\varphi)\in\big(0,\frac{-1+\sqrt{1+\lambda^{2}}}{\lambda})$.

\begin{lemma}\label{peri-lem}
For a reduced spherical polygon $V=v_{1}v_{2}\cdots v_{n}$ with thickness less than $\frac{\pi}{2}$,
the perimeter is $\mathrm{perim}(V)=2\sum _{i=1}^{n}f_{1}(y_{i})$, where $y_{i}=g(\varphi_{i})$.
\end{lemma}
\proof
For each $i\in\{1,2,\dots,n\}$, we focus on the right spherical triangle $o_{i}t_{i}v_{i+(n+1)/2}$ in $V$.
Let $|t_{i}v_{i+(n+1)/2}|=a_{i}$, $|o_{i}t_{i}|=b_{i}$ and $|o_{i}v_{i+(n+1)/2}|=c_{i}$.
By Lemma \ref{cog-tri}, we obtain
$\mathrm{perim}(V)=2\sum _{i=1}^{n}|v_{i}t_{i+(n+1)/2}|=2\sum _{i=1}^{n} a_{i}$.
Moreover, we have $|o_{i}t_{i}|+|o_{i}v_{i+(n+1)/2}|=b_{i}+c_{i}=\omega$.
Note that $b_{i}<\omega$.
From Equation \eqref{eq1}, we obtain
\begin{equation*}
\cos\varphi_{i}=\frac{\tan b_{i}}{\tan c_{i}}=\frac{\tan b_{i}}{\tan(\omega-b_{i})}
=\frac{\tan b_{i}(1+\tan\omega\tan b_{i})}{\tan\omega-\tan b_{i}}.
\end{equation*}
Let $\tan b_{i}=y_{i}$, from this and $\lambda=\tan\omega$, the above equation becomes
\begin{equation}\label{in1}
\cos\varphi_{i}=\frac{y_{i}(1+\lambda y_{i})}{\lambda-y_{i}},
\end{equation}
thus
\begin{equation}\label{lem-inq2}
 \varphi_{i}=\arccos\frac{y_{i}(1+\lambda y_{i})}{\lambda-y_{i}}.
\end{equation}
Moreover, we have
$$
\tan c_{i}=\frac{\tan b_{i}}{\cos\varphi_{i}}=\frac{\lambda-y_{i}}{1+\lambda y_{i}},
$$
and thus $c_{i}=\arctan\frac{\lambda-y_{i}}{1+\lambda y_{i}}$.
By Equation \eqref{eq2}, we get that
$$
\cos a_{i}=\frac{\cos c_{i}}{\cos b_{i}}
=\frac{1+\lambda y_{i}}{\sqrt{1+\lambda^{2}}},
$$
hence $a_{i}=\arccos\frac{1+\lambda y_{i}}{\sqrt{1+\lambda^{2}}}=f_{1}(y_{i})$.
Therefore, $\mathrm{perim}(V)=2\sum _{i=1}^{n}f_{1}(y_{i})$,
where $$y_{i}=\frac{-(1+\cos\varphi_{i})+\sqrt{(1+\cos\varphi_{i})^{2}+4\lambda^{2}\cos \varphi_{i}}}{2\lambda}=g(\varphi_{i}),$$ which can be calculated  by Equation \eqref{in1}.
\hfill\qed

\begin{theorem}\label{thm-perim}
The regular spherical $n$-gon has the minimum perimeter among all regular spherical $k$-gons of fixed thickness $\omega\in(0,\frac{\pi}{2})$, where $k,n$ are odd and $3\leq k\leq n$.
\end{theorem}
\proof
Let $V_{k}=v_{1}v_{2}\cdots v_{k}$ be a regular spherical odd-gon.
Lemma \ref{col3.3} shows that $V_{k}$ is reduced, then by Lemma \ref{peri-lem}, we have
$$\mathrm{perim}(V_{k})=2\sum _{i=1}^{k}f_{1}(y_{i}),$$ where $y_{i}=g(\varphi_{i})$.
$V_{k}$ is a regular polygon, then Lemma \ref{lem3.5} $(2)$ shows that $\varphi_{1}=\cdots=\varphi_{k}=\frac{\pi}{k}$.
Thus we have $y_{1}=\cdots=y_{k}=g(\frac{\pi}{k})$, where
$$
g(\frac{\pi}{k})=\frac{-(1+\cos\frac{\pi}{k})+\sqrt{(1+\cos\frac{\pi}{k})^{2}+4\lambda^{2}\cos \frac{\pi}{k}}}{2\lambda}.
$$
For brevity, we denote $\varphi_{i}$ and $y_{i}$ by $\varphi$ and $y$, respectively, where $i=1,2,\ldots,k$.
Then $\varphi=\frac{\pi}{k}$ and $y=g(\varphi)$, where $y\in(0,\frac{-1+\sqrt{1+\lambda^{2}}}{\lambda})$.
Equation \eqref{lem-inq2} shows that $$\varphi=\arccos\frac{y(1+\lambda y)}{\lambda-y}=f_{2}(y),$$
then $k=\frac{\pi}{\varphi}=\frac{\pi}{f_{2}(y)}$.
As a result,
\begin{equation}\label{reg-perim}
\mathrm{perim}(V_{k})=2kf_{1}(y)=2\pi\frac{f_{1}(y)}{f_{2}(y)},
\end{equation}
where $y\in(0,\frac{-1+\sqrt{1+\lambda^{2}}}{\lambda})$.

Lemma \ref{dec-func} permits that $\frac{f_{1}(y)}{f_{2}(y)}$ is a decreasing function of $y$.
Since $$y=g(\frac{\pi}{k})=\frac{-(1+\cos\frac{\pi}{k})+\sqrt{(1+\cos\frac{\pi}{k})^{2}+4\lambda^{2}\cos \frac{\pi}{k}}}{2\lambda},$$
it follows that $y$ is an increasing function of $k$.
Thus the perimeter of the regular spherical polygon decreases with the increase of $k$. The conclusion is as desired.
\hfill\qed

\begin{theorem}\label{regular-min-perim}
The regular spherical $n$-gon has the minimum perimeter among all
reduced spherical $n$-gons with the same thickness $\omega\in(0,\frac{\pi}{2})$.
\end{theorem}

\proof
Let $V=v_{1}v_{2}\cdots v_{n}$ be an arbitrary reduced spherical $n$-gon,
then by Lemma \ref{peri-lem}, the perimeter of $V$ is $\mathrm{perim}(V)=2\sum _{i=1}^{n}f_{1}(y_{i})$, where $y_{i}=g(\varphi_{i})$ and
$$g(\varphi_{i})=\frac{-(1+\cos\varphi_{i})+\sqrt{(1+\cos\varphi_{i})^{2}+4\lambda^{2}\cos \varphi_{i}}}{2\lambda}.$$
Let $F(\varphi_{i})=f_{1}(g(\varphi_{i}))=\arccos\frac{1+\lambda g(\varphi_{i})}{\sqrt{1+\lambda^{2}}}$, where $\varphi_{i}\in(0,\frac{\pi}{2})$ and $\lambda=\tan\omega\in(0,+\infty)$.
Thus $\mathrm{perim}(V)=2\sum _{i=1}^{n}F(\varphi_{i})$.
Lemma \ref{F(z)} states that $F(x)$ is a convex function of $x$.
Thus from Jensen's inequality \cite{Jensen}, we obtain that
$$\frac{F(\varphi_{1})+\cdots+F(\varphi_{n})}{n}
\geq F(\frac{\varphi_{1}+\cdots +\varphi_{n}}{n}),$$
the equality holds when $\varphi_{1}=\cdots=\varphi_{n}$.
Then the perimeter of $V$ satisfies
$$\mathrm{perim}(V)=2n(\frac{F(\varphi_{1})+\cdots+F(\varphi_{n})}{n})
\geq 2nF(\frac{\varphi_{1}+\cdots +\varphi_{n}}{n}).$$

{\bf Case 1.} In the case $V$ is a regular spherical polygon. By Equation \eqref{reg-perim}, it follows that $\mathrm{perim}(V)=2nF(\frac{\pi}{n})$.

{\bf Case 2.} In the case $V$ is a non-regular spherical polygon.
Lemmas \ref{fact-2020} and \ref{lem3.5} imply that
$$\pi\leq\varphi_{1}+\cdots+\varphi_{n}<\frac{n\pi}{2}.$$
Since $F(x)$ is an increasing function of $x$ within the domain $x\in(0,\frac{\pi}{2})$,
it follows that $$F(\frac{\varphi_{1}+\cdots +\varphi_{n}}{n})\geq F(\frac{\pi}{n}).$$
Consequently, we get
$$\mathrm{perim}(V)=2n(\frac{F(\varphi_{1})+\cdots+F(\varphi_{n})}{n})
\geq2nF(\frac{\varphi_{1}+\cdots +\varphi_{n}}{n})\geq 2nF(\frac{\pi}{n}).$$

The conclusion follows by the above two cases.
\hfill\qed

By Theorems \ref{thm-perim} and \ref{regular-min-perim}, we confirm the Problem $\mathrm{P1}$.

\section{The maximal diameter of reduced spherical polygons}
Recall that the Problem $\mathrm{P2}$ is proposed by Lassak after the discussion of
\cite[Theorem $4.2$]{polygon}. In the following, we establish some related conclusions. Denote by $\mathrm{diam}(V)$ the diameter of a spherical polygon $V$.

\begin{lemma}{\rm\cite[Proposition 4.1]{polygon}}\label{prop4.1}
The diameter of any reduced spherical $n$-gon is realized only for some pairs of vertices whose indices $\pmod {n}$ differ by $\frac{n-1}{2}$ or $\frac{n+1}{2}$.
\end{lemma}

\begin{lemma}{\rm\cite[Theorem 4.2]{polygon}}\label{Thm4.2}
For every reduced spherical polygon $V$ on the sphere, we have
\begin{equation}\label{inq-dim}
 \mathrm{diam}(V)\leq \arccos\Big(\cos\Delta(V)\sqrt{1-\frac{\sqrt{2}}{2}\sin\Delta(V)}~\Big),
\end{equation}
with equality for the regular spherical triangle in the part of $V$.
\end{lemma}

The Problem $\mathrm{P2}$ asks if the equality in Lemma \ref{Thm4.2} holds only for regular triangle? The answer is negative, so we derive a new inequality different from \eqref{inq-dim} and confirm the Problem $\mathrm{P2}$ within our investigation.
The following statement is analogous to the second part of \cite[Theorem $9$]{1990}.
Although the proof process is trivial and the method is the same as \cite[Theorem $4.2$]{polygon}, for the completeness, we restate as follows.

\begin{theorem}\label{diameter}
For every reduced spherical polygon $V=v_{1}\cdots v_{n}$ of thickness less than $\frac{\pi}{2}$, we have
\begin{equation}\label{inq-dim-new}
\mathrm{diam}(V)\leq 2\arccos(\frac{1}{2\sin\gamma}),
\end{equation}
where $\gamma=\arcsin\frac{-\cos\Delta(V)+\sqrt{\cos^{2}\Delta(V)+8}}{4}$,
and the equality holds only for the regular spherical triangle.

\end{theorem}
\proof
By Lemma \ref{prop4.1}, we assume that $\mathrm{diam}(V)=|v_{i}v_{i+(n+1)/2}|$ for some $i\in\{1,2,\ldots,n\}$ (in the case $\mathrm{diam}(V)=|v_{i}v_{i+(n-1)/2}|$, the following is similar).

Consider the spherical triangle $v_{i}t_{i}v_{i+(n+1)/2}$.
Put
$$
a_{i}=|t_{i}v_{i+(n+1)/2}|,~c_{i}=|v_{i}v_{i+(n+1)/2}|,~\mathrm{and}~
\gamma_{i}=\angle v_{i}v_{i+(n+1)/2}t_{i}.$$
Lemmas \ref{col3.9} and \ref{cog-tri} (resp. \cite[Corollary $3.7$]{polygon}) imply that $\gamma_{i}=\alpha_{i}+\beta_{i}\geq 2\beta_{i}$, where $0<\gamma_{i}\leq\frac{\pi}{2}$ (for the range of $\gamma_{i}$, we can refer to the proof process of Lemma \ref{relation}).
From $|v_{i}t_{i}|=\Delta(V)$ and Equation \eqref{eq-4}, we obtain
$$\frac{\sin\beta_{i}}{\sin a_{i}}=\frac{\sin\gamma_{i}}{\sin \Delta(V)}.$$
Thus from $\gamma_{i}\geq 2\beta_{i}$, we acquire
$\frac{\sin\beta_{i}}{\sin a_{i}}\geq\frac{\sin 2\beta_{i}}{\sin\Delta(V)}$
and hence
\begin{equation}\label{diam-inq1}
 \sin a_{i}\leq \frac{\sin\Delta(V)}{2\cos\beta_{i}}.
\end{equation}
Moreover, by Lemma \ref{relation}, we have $\beta_{i}\leq\gamma$ (In the proof process of Lemma \ref{Thm4.2}, Lassak use a slight rough upper bound $\beta_{i}\leq\frac{\pi}{4}$ to obtain the final inequality. This is the main difference between our proof with Lassak's).
Hence we get $\sin a_{i}\leq\frac{\sin\Delta(V)}{2\cos\gamma}$
and then
\begin{equation}\label{diam-inq2}
\cos a_{i}\geq\sqrt{1-\frac{\sin^{2}\Delta(V)}{4\cos^{2}\gamma}}.
\end{equation}
Applying the Pythagorean theorem (i.e. Equation \eqref{eq2}) on sphere to \eqref{diam-inq2},  we obtain
\begin{equation*}
\cos c_{i}\geq\cos\Delta(V)\sqrt{1-\frac{\sin^{2}\Delta(V)}{4\cos^{2}\gamma}},
\end{equation*}
as a result,
\begin{equation*}
\mathrm{diam}(V)=c_{i}\leq \arccos\Big(\cos\Delta(V)\sqrt{1-\frac{\sin^{2}\Delta(V)}{4\cos^{2}\gamma}}~\Big).
\end{equation*}

Observe that the equality holds only when \eqref{diam-inq1} and \eqref{diam-inq2} attain the equality, which requires $\gamma_{i}=2\beta_{i}$ and $\beta_{i}=\gamma$. This leads to the equality $\gamma=\beta_{i}=\alpha_{i}$ and thus $V$ is a regular spherical triangle by Lemma \ref{relation}.
From the above discussions and Lemma \ref{gamma-regular triangle} $(2)$, we obtain the following equality
$$\arccos\Big(\cos\Delta(V)\sqrt{1-\frac{\sin^{2}\Delta(V)}{4\cos^{2}\gamma}}~\Big)=2\arccos(\frac{1}{2\sin\gamma}).$$
This completes the proof.
\hfill\qed

Lemma \ref{col3.3} implies that the thickness of every reduced spherical polygon is at most $\frac{\pi}{2}$. Then in Lemma \ref{Thm4.2}, the thickness of $V$ satisfies $\Delta(V)\leq\frac{\pi}{2}$, while there exists an error in this lemma. For clarity, we establish the following version.

\begin{proposition}\label{Thm4.2-modify}
For every reduced spherical polygon $V$ of thickness at most $\frac{\pi}{2}$, we have
\begin{equation}\label{inq-dim-modify}
 \mathrm{diam}(V)\leq \arccos\Big(\cos\Delta(V)\sqrt{1-\frac{\sqrt{2}}{2}\sin\Delta(V)}~\Big),
\end{equation}
the equality holds only when the thickness of $V$ is $\frac{\pi}{2}$.
\end{proposition}
\proof
Since \eqref{inq-dim-modify} has been proved in \cite[Theorem $4.2$]{polygon}, we omit it.

In the case $\Delta(V)<\frac{\pi}{2}$, by Theorem \ref{diameter}, it follows that $\mathrm{diam}(V)\leq 2\arccos(\frac{1}{2\sin\gamma})$. One can check that
$$2\arccos(\frac{1}{2\sin\gamma})< \arccos\Big(\cos\Delta(V)\sqrt{1-\frac{\sqrt{2}}{2}\sin\Delta(V)}~\Big).$$

In the case $\Delta(V)=\frac{\pi}{2}$, \cite[Proposition $2$]{Lassak2020} implies that $V$ is a body of constant width, thus $\mathrm{diam}(V)=\Delta(V)=\frac{\pi}{2}$.
One can also calculate that the right hand side of the inequality \eqref{inq-dim-modify} equals $\frac{\pi}{2}$.
\hfill\qed

By the above investigation, Problem $\mathrm{P2}$ holds only when considering Theorem \ref{diameter}.

\section{The smallest radius of a spherical disk covering reduced spherical polygons}
The problems about covering the reduced convex bodies by a disk in $E^{2}$ and $S^{2}$ can be found in \cite{2003} and \cite{Michal}, respectively.
For the planar case \cite{2003}, every reduced body $R$ (resp. reduced convex polygon) is contained in a disk of radius $\frac{\sqrt{2}}{2}\Delta(R)$ (resp. $\frac{2}{3}\Delta(R)$).

In the sphere, \cite[Theorem $2$]{Michal} shows that every reduced spherical body $R$ of thickness at most $\frac{\pi}{2}$ is contained in a disk of radius $\arctan(\sqrt{2}\tan\frac{\Delta(R)}{2})$.
And we can not improve the above estimate when considering all reduced spherical bodies.
Now we concentrate on the Problem $\mathrm{P3}$ to find the smallest radius of a disk that contains every reduced polygon of a given thickness on $S^{2}$.

In \cite{1995Jung}, Dekster extends the Jung Theorem to the $n$-dimensional hyperbolic space and the $n$-dimensional sphere ($n\geq2$). We only exhibit the conclusions we need.

\begin{lemma}{\rm\cite[Theorem 2]{1995Jung}}\label{Jung-Thm}
Let $S$ be a compact set in $S^{n}$ of diameter $D$ and circumradius $r$. Let $B$ be a metric ball of radius $r$ containing $S$. Then
$$
D\geq 2\arcsin(\sqrt{\frac{n+1}{2n}}\sin r),~r\in[0,\pi],~for~S\subseteq S^{n}.
$$
\end{lemma}

This version applies to compact sets. The reduced spherically convex polygon is exactly a compact (bounded and closed) set in $S^{2}$. So we could get the following statement.

\begin{theorem}\label{cover-disk}
Every reduced spherically convex polygon $V$ of thickness $\Delta(V)\in(0,\frac{\pi}{2})$ is contained in a disk of radius $\arcsin(\frac{2}{\sqrt{3}}\sqrt{1-\frac{1}{4\sin^{2}\gamma}}~)$, where
$$ \gamma=\arcsin\frac{-\cos\Delta(V)+\sqrt{\cos^{2}\Delta(V)+8}}{4}.
$$
\end{theorem}
\proof Denote by $r$ the circumradius of $V$. In $S^{2}$, by Lemma \ref{Jung-Thm}, we have $\mathrm{diam}(V)\geq 2\arcsin(\frac{\sqrt{3}}{2}\sin r)$,
thus
$$r\leq\arcsin(\frac{2}{\sqrt{3}}\sin\frac{\mathrm{diam}(V)}{2}).$$
Combining this with \eqref{inq-dim-new} of Theorem \ref{diameter}, it follows that $$r\leq\arcsin(\frac{2}{\sqrt{3}}\sin\frac{2\arccos(\frac{1}{2\sin\gamma})}{2})=\arcsin\Big(\frac{2}{\sqrt{3}}\sqrt{1-\frac{1}{4\sin^{2}\gamma}}~\Big).$$
The conclusion is as desired.
\hfill\qed

Note that the inequality \eqref{inq-dim-new} of Theorem \ref{diameter} attains equality only when $V$ is a regular triangle, so the estimate $\arcsin(\frac{2}{\sqrt{3}}\sqrt{1-\frac{1}{4\sin^{2}\gamma}})$ is sharp only for the regular triangle, and the estimate can not be improved.

The following table lists some values of the radius under specific thickness.

\begin{table}[h]
\centering
\resizebox{\textwidth}{!}{%
\begin{tabular}{|c|c|c|c|c|}
\hline
radius $\backslash$ thickness & $\frac{\pi}{8}$            & $\frac{\pi}{6}$     & $\frac{\pi}{4}$       & $\frac{\pi}{3}$      \\ \hline
$\arcsin(\frac{2}{\sqrt{3}}\sqrt{1-\frac{1}{4\sin^{2}\gamma}})$& $0.260304\cdots$& $0.345523\cdots$ &$ 0.511669\cdots$ & $0.670020\cdots$ \\ \hline
\end{tabular}}
\caption{The values of radius under specific thickness}
\label{table}
\end{table}

\end{document}